\baselineskip=15pt plus 2pt
\magnification=1150
\def\sqr#1#2{{\vcenter{\vbox{\hrule height.#2pt\hbox{\vrule width.#2pt
height#1pt\kern#1pt \vrule width.#2pt}\hrule height.#2pt}}}}
\def\square{\mathchoice\sqr64\sqr64\sqr{2.1}3\sqr{1.5}3}
\font\medtenrm=cmr10 scaled \magstep2

\centerline {\medtenrm A logarithmic Sobolev inequality for the
invariant measure}\par
\centerline {\medtenrm of the periodic Korteweg--de Vries equation}\par
\vskip.05in
\centerline {\medtenrm Gordon Blower}\par
\centerline {\sl Department of Mathematics and Statistics}\par
\centerline {\sl Lancaster University}\par
\centerline {\sl Lancaster LA1 4YF}\par
\centerline {\sl England}\par
\centerline {Email g.blower@lancaster.ac.uk}\par
\vskip.1in
\centerline {21st December 2009}\par
\vskip.1in
\hrule
\vskip.1in

\noindent {\bf Abstract}\par
\noindent  The periodic KdV equation
$u_t=u_{xxx}+\beta uu_x$
arises from a Hamiltonian system with infinite-dimensional phase space
$L^2({\bf T})$.
Bourgain has shown that there exists a Gibbs probability measure $\nu$
on balls $\{ \phi :\Vert \phi\Vert^2_{L^2}\leq N\}$ in the phase space
 such that the Cauchy problem for KdV is well posed
on the support of $\nu$, and $\nu$ is invariant under the KdV
flow. This paper shows that $\nu$ satisfies a logarithmic Sobolev 
inequality. The stationary points of the Hamiltonian on spheres are
found in terms of elliptic functions, and they are shown to be
linearly stable. The paper also presents logarithmic Sobolev
inequalities for the modified periodic KdV equation and the cubic
nonlinear Schr\"odinger equation, for small values of $N$.\par
\vskip.1in
\noindent {\bf R\'esum\'e}\par
\noindent L'\'equation KdV p\'eriodique $u_t=u_{xxx}+\beta uu_x$ 
r\'esulte d'un syst\`eme
hamiltonien des espaces de phases  $L^2({\bf T})$ de dimension
infinite. Bourgain a d\'emontr\'e qu'il existe une mesure de probabilit\'e de Gibbs
$\nu$ sur les boules  $\{ \phi :\Vert \phi\Vert^2_{L^2}\leq N\}$ dans 
l'espace des phases telles que le probl\`eme
de Cauchy pour KdV est bein pos\'e sur le support de $\nu$ et $\nu$ est
invariante sous le flux de KdV. Cet article montre que $\nu$ satisfait
 une in\'egalit\'e de Sobolev logarithmique. Les points fixes de
l'hamiltonien sur les sph\`eres sont trouv\'es en termes de fonctions
elliptiques, et il est d\'emontr\'e qu'elles sont lin\'eairement 
stables. L'article pr\'esente \'egalement les in\'egalit\'es de Sobolev
logarithmique pour l'\'equation KdV modifie et les cubiques
non lin\'eares
de Schr\"odinger pour les petites valeurs de $N$.

\par
\vskip.1in
\hrule
\vskip.1in
\noindent {\bf 1. Introduction}\par
\vskip.05in
\noindent In this paper, we are concerned with solutions of the KdV
equation that are periodic in the space variable and typical in the
sense that they form the support of an invariant measure on an
infinite-dimensional phase space. Specifically, we consider $u:{\bf
T}\times (0, \infty )\rightarrow {\bf R}$ such that $u(\, ,t)\in
L^2({\bf T})$ for each $t>0$, then we introduce the Hamiltonian
$$H(u)={{1}\over{2}}\int_{\bf T}\Bigl({{\partial u}\over{\partial 
x}}(x,
t)\Bigr)^2{{dx}\over{2\pi}}-{{\beta}\over{6}}\int_{\bf T}
u(x, t)^3 {{dx}\over{2\pi}}.\eqno(1.1)$$
\noindent Here $\beta$ is the reciprocal of temperature, and without
loss of generality, we assume throughout that $\beta >0$. The canonical equation of motion is
$${{\partial u}\over{\partial
t}}={{\partial}\over{\partial x}}{{\delta H}\over{\delta
u}},\eqno(1.2)$$
\noindent which gives the Korteweg--de Vries equation
$${{\partial u}\over{\partial t}}=-{{\partial^3u}\over{\partial
x^3}}-\beta u{{\partial u}\over{\partial x}}.\eqno(1.3)$$
\noindent Given a solution of (1.3) that is suitably differentiable, 
one can easily verify that\par
\noindent $\int_{\bf T}
u(x,t)^2dx/(2\pi)$ and $H(u)$ are invariant with respect
to time. In order to ensure that the Gibbs measure can be
normalized, we work on bounded subsets of $L^2({\bf T})$. As in [4], we 
introduce the particle number $N<\infty$, the ball 
$$B_N=\Bigl\{\phi \in L^2({\bf T}): \int_{\bf T}\phi (x
)^2{{dx}\over{2\pi}}\leq N\Bigr\}\eqno(1.4)$$
\noindent with indicator function ${\bf I}_{B_N}$ and the Gibbs measure 
$$\nu_N^\beta (d\phi )=Z_N(\beta )^{-1}{\bf I}_{B_N}(\phi)e^{-H(\phi
)}\prod_{e^{ix}\in {\bf T}} d\phi (x )\eqno(1.5)$$
\noindent where the normalizing constant $Z_N(\beta )$ is so chosen as to give
a probability measure.\par
\vskip.05in
\noindent {\bf Definition.} The modified canonical ensemble is the 
probability space $(B_N, \nu_N^\beta )$
that has particle number $N$ at
inverse temperature $\beta$.\par
\vskip.05in
\indent In view of [12], the natural canonical ensemble would appear
to be a probability measure on
the sphere $S_N=\{
\phi\in L^2({\bf T}): \Vert \phi\Vert_{L^2}^2=N\}$, but this is
technically difficult to deal with, so we prefer the modified canonical
ensemble. However, in section 3 we consider the Hamiltonian on the
sphere and show that the stationary points of $H$ on $S_N$ are
given by elliptic functions.\par  
\indent There are various means for introducing Gibbs measures on
infinite-dimensional phase spaces. In [8], Lebowitz, Rose and Speer constructed an invariant
measure for the nonlinear Schr\"odinger equation on the line, and
investigated the stability of the ground state. Using purely
probabilistic arguments, McKean and Vaninsky gave an alternative
construction [14].\par
\indent Here we construct the measure 
via random Fourier series. We write $\phi (x)\sim
{{1}\over{2}}a_0+\sum_{j=1}^\infty (a_j\cos jx +b_j\sin jx
)$, and regard $(a_j, b_j)$ as an $\ell^2$ sequence of coordinates for
$\phi\in L^2({\bf T})$. Let $(\gamma_j)_{j=-\infty}^\infty$ be mutually
independent standard Gaussian random variables on some
probability space $(\Omega , {\bf P})$, and let $W$ be the probability
measure on $L^2({\bf T})$ that is induced by 
$$\omega\mapsto \phi_\omega (x )=\gamma_0+\sum_{j=-\infty}^{-1} \gamma_j
{{\sin jx }\over{j}}+\sum_{j=1}^\infty \gamma_j {{\cos
jx}\over{j}},\eqno(1.6)$$
\noindent namely Brownian loop. Then $\Omega_N=\{ \omega\in \Omega:
{{1}\over{2}}\sum_{j=-\infty ;j\neq 0}^\infty \gamma_j^2/j^2\leq N\}$
maps into $B_N$, and we can introduce the Gibbs measure as
$$\nu_N^\beta (d\phi )=Z_N(\beta )^{-1}{\bf I}_{\Omega_N}(\omega )\exp
\Bigl( {{\beta}\over{6}}\int_{\bf T}\phi_\omega (x
)^3{{dx}\over{2\pi}}\Bigr) W(d\phi_\omega ).\eqno(1.7)$$
\noindent Bourgain [2, 3] shows that there exists $Z_N(\beta )>0$ such that
$\nu_N^\beta $ is a Radon probability measure on the closed subset $B_N$ of $L^2({\bf
T})$. Further, the Cauchy initial value problem
$$\cases{ u_t=-u_{xxx}-\beta uu_x\cr
          u(x,0)=\phi (x)\cr}\eqno(1.8)$$
\noindent is locally well posed on the support of $\nu_N^\beta$; more
precisely, for each $\delta>0$, there exists $\tau (\delta )>0$ and a
compact set $K_\delta$ such that $\nu_N^\beta (K_\delta )>1-\delta $ and
such that for all
$\phi\in K_\delta$ there exists a unique
solution $u(x,t)$ to (1.8) for $t\in [0, \tau (\delta )]$. Existence
of the invariant measure $\nu_N^\beta$ implies that the local
solution extends to a global solution for almost all initial data
with respect to $\nu_N^\beta$. We should expect the long term
behaviour of solutions to
consist of a solitary travelling wave coupled with fluctuations, as
described by the invariant measure. The main result of this paper is a logarithmic Sobolev
inequality which shows that such a space of solutions is stable.\par
\par
\vskip.05in
\noindent {\bf Definition.} Suppose that $F:B_N\rightarrow {\bf R}$ is 
G\^ateaux
differentiable, so that for all $\phi$ inside $B_N$, there exists 
$\nabla F(\phi )\in L^2({\bf T})$ such that 
$$\langle \nabla F(\phi ), \psi \rangle_{L^2}=\lim_{t\rightarrow 0+}
{{F(\phi +t\psi )-F(\phi )}\over{t}}\eqno(1.9)$$
\noindent for all $\psi\in L^2$. Suppose further that the limit
exists uniformly on $\{ \psi\in L^2: \Vert\psi\Vert_{L^2}=1\}$; then
$F$ is Fr\'echet differentiable.\par
\vskip.05in
\indent Let $\dot H^{1/2}=\{ \phi (x)=\sum_{k\neq 0; k=-\infty}^\infty
a_ke^{ikx}:\sum_{k\neq 0; k=-\infty }^\infty \vert k\vert \vert
a_k\vert^2<\infty \}$, and let $G:L^2({\bf T})\rightarrow L^2({\bf T})$
be the operator
$$G\phi (x)=\int_{\bf T}\log {{1}\over{\vert e^{ix}-e^{iy}\vert}}\,
\phi (y)\, {{dy}\over{2\pi}}\sim \sum_{k\neq 0;k=-\infty}^\infty 
{{\hat \phi (k)}\over{\vert k\vert}} e^{ikx}.\eqno(1.10)$$
\noindent Then
$$\langle G\psi , \phi \rangle_{\dot H^{1/2}}=\langle \psi
,\phi\rangle_{L^2}.\eqno(1.11)$$
\noindent We write $\delta F(\phi )=G(\nabla F(\phi )),$ and observe
that $\Vert\delta F(\phi )\Vert_{L^2}\leq \Vert\nabla F(\phi
)\Vert_{L^2}.$\par

\vskip.05in
\noindent {\bf Definition} {(\sl Logarithmic Sobolev inequality). [16]} Say
that a probability measure $\nu_N$ on $B_N$ satisfies the logarithmic
Sobolev inequality with constant $\alpha >0$ if  
$$\int_{B_N}F(\phi )^2\log \Bigl(F(\phi
)^2/\int_{B_N}F^2d\nu_N\Bigr) \nu_N^\beta (d\phi )\leq
{{2}\over{\alpha}}\int_{B_N}\Vert\delta F(\phi )\Vert^2_{L^2}
\nu_N(d\phi
)\eqno(1.12)$$
\noindent for all Fr\'echet differentiable functions $F\in L^2(B_N; \nu_N^\beta )$ such
that 
$\Vert \delta F\Vert_{L^2}\in L^2(B_N;\nu_N
)$.\par

\noindent {\bf Theorem 1.} {\sl For all $\beta ,N>0$ the measure
 $\nu_N^\beta$ satisfies the logarithmic Sobolev inequality with 
$$\alpha =2^{-1}\exp\bigl( -C\beta^{5/2}N^{9/4}\bigr)\eqno(1.13)$$
some absolute 
constant $C>0$.}\par

\vskip.05in
\indent In section 2 we prove Theorem 1 and deduce a concentration
inequality concerning Lipschitz functions on the $B_N$. This shows that
certain random variables are tightly concentrated around their mean
values, just as a Gaussian random variable is concentrated close to its
mean. In
section 3, we consider stability of the stationary points of the
Hamiltonian restricted to spheres. Analysis of the stationary points
reduces to classical spectral theory of Lam\'e's equation, and we are able to identify
stationary points as elliptic functions.\par 
\indent In section 4, we consider the Gibbs measure associated with the
modified periodic KdV equation, and obtain a logarithmic Sobolev
inequality when $N\beta$ is small and positive. Finally, in section 5
we apply a similar
analysis to the periodic cubic Schr\"odinger equation.\par  
\vskip.1in
\noindent {\bf 2. Logarithmic Sobolev inequalities}\par
\vskip.05in
\noindent In order to prove Theorem 1, we express the phase space in
terms of Fourier components. As in Parseval's identity, there is a unitary map
$\ell^2\rightarrow L^2({\bf T})$ 
$$(a_0; a_n, b_n)_{n=1}^\infty \mapsto a_0+\sum_{k=1}^\infty \sqrt{2}( 
a_k\cos kx +b_k\sin kx);\eqno(2.1)$$
\noindent under this correspondence, $F:B_N\rightarrow {\bf R}$ may be
identified with $f:\Omega_N\rightarrow {\bf R}$ and $\nabla F$ corresponds
to $( {{\partial f}\over{\partial a_j}}, {{\partial f}\over{\partial
b_{j+1}}})_{j=0}^\infty.$ 
To see this, we consider $\psi (x)=
\sum_{k=1}^\infty \sqrt{2}( 
c_k\cos kx +d_k\sin kx)$ and observe that 
$$\langle \nabla f,
(c_k,d_k)_{k=1}^\infty\rangle_{\ell^2}=\sum_{k=1}^\infty\Bigl( {{\partial
f}\over{\partial a_k}}c_k+{{\partial f}\over{\partial
b_k}}d_k\Bigr)\eqno(2.2)$$
\noindent while 
$$\langle \nabla F, \psi\rangle_{L^2}=\int (\nabla F)\bar \psi 
{{dx}\over{2\pi}},\eqno(2.3)$$
\noindent and we can recover the Fourier coefficients of
$\nabla F$. Further, the Gibbs measure $\nu_N^\beta$ may
be expressed in terms of the Fourier components as    
$$Z_N(\beta )^{-1} \exp\Bigl[ {{\beta}\over{6}}\int_{\bf T}\Bigl(
a_0+\sqrt{2}\sum_{j=1}^\infty (a_j\cos jx +b_j\sin jx
)\Bigr)^3{{dx}\over{2\pi}}-a_0^2-\sum_{j=1}^\infty j^2(a_j^2+b_j^2)
\Bigr]$$
$$\times {\bf I}_{[0, N]}\Bigl( a_0^2+\sum_{j=1}^\infty (a_j^2+b_j^2)
\Bigr){{da_0}\over{\sqrt{2\pi}}}\prod_{j=1}^\infty
{{j^2 da_jdb_j}\over{2\pi}}.\eqno(2.4)$$
\noindent For notational simplicity, we write $a_{-j}=b_j$ for
$j=1, 2,\dots $, and assume that $0<\beta <\sqrt {3}/(4\pi\sqrt{N})$.
We introduce the potential
$$V (a, b)={{a_0^2}\over{2}}+{{1}\over{2}}
\sum_{j=1}^\infty j^2(a_j^2+b_j^2)-{{\beta
2\sqrt{2}}\over{6}}\int_{\bf T}\Bigl({{a_0}\over{\sqrt{2}}}+ \sum_{j=1}^\infty 
(a_j\cos jx +b_j\sin jx )\Bigr)^3{{dx}\over {2\pi}}.\eqno(2.5)$$
\vskip.05in
\noindent {\bf Lemma 2.} {\sl Suppose that $0<\beta \sqrt N\leq
\sqrt{3} /(32\pi )$. Then $V$ is uniformly convex on $B_N$ and has a unique minimum
at the origin; moreover, $\nu_N^\beta$ satisfies the logarithmic 
Sobolev inequality (1.12) with $\alpha =1/2$.}

\vskip.05in

\noindent {\bf Proof.} First we scale the variables to $x_j=ja_j$ and $y_j=jb_j$, so
that 
$$\Omega_N=\Bigl\{ (a_0;a_j,b_j)\in \ell^2: a_0^2+\sum_{j=1}^\infty
(a_j^2+b_j^2)\leq N\Bigr\}\eqno(2.6)$$ 
\noindent is transformed to the ellipsoid 
$${\cal E}_N=\Bigl\{ (x_0; x_j,
y_j): x_0^2+\sum_{j=1}^\infty (x_j^2+y_j^2)/j^2\leq N\Bigr\}.$$
\noindent Let $G:{\cal E}_N\rightarrow \Omega_N$ be the diagonal map
$G:((x_k,y_k))_{k=1}^\infty\mapsto ((x_k/k, y_k/k))_{k=1}^\infty$, with
left inverse $D:((a_k,b_k))_{k=1}^\infty =((ka_k, kb_k))_{k=1}^\infty$ so
that $DG=I$. We then introduce $W:{\cal E}_N\rightarrow {\bf R}$ by
$W(x)=V(G(x))$. \par
\indent To verify Bakry and Emery's criterion [1] for the logarithmic
Sobolev inequality, we need to show
that the Hessian matrix of $W$ satisfies
$${\hbox{Hess}}\, W\geq
{{1}\over{2}}I\eqno(2.7)$$
\noindent  and hence that
$\omega_N^\beta =\zeta_N(\beta )^{-1}\exp (-W(x))\, dx$ satisfies the
logarithmic Sobolev inequality
$$\int_{{\cal E}_N} g(x)^2\log \Bigl( g(x)^2/\int g^2d\omega_N^\beta
\Bigr)\, dx\leq 4
\int_{{\cal E}_N} \Vert \nabla g\Vert_{\ell^2}^2 d\omega_N^\beta .\eqno(2.8)$$
\noindent Now $G$ induces $\nu_N^\beta $ from $\omega_N^\beta$; so with
$g=f\circ G$ we have $\nabla g=((\nabla f)\circ G)(\nabla G )$ where
$\nabla G$ is represented by the diagonal matrix $(1/k)_{k=1}^\infty$.
The condition (2.7) is equivalent to 
$${\hbox {Hess}}\, V=\left[\matrix{ 
{{\partial^2V}\over{\partial a_j\partial
a_k}}&{{\partial^2V}\over{\partial a_j\partial b_k}}\cr
{{\partial^2V}\over{\partial a_j\partial
b_k}}&{{\partial^2V}\over{\partial b_j\partial b_k}}\cr}\right]\geq
{{1}\over{2}}D^2\eqno(2.9)$$
Let $D$ be the diagonal matrix $(j)$ with respect to the
Fourier basis and let
$$\eqalignno{ v_{jk}&={{\partial^2}\over{\partial a_j\partial a_k}}{{2\sqrt
2}\over{6}}\int_{\bf T} \Bigl({{a_0}\over{\sqrt{2}}}+ \sum_{\ell =1}^\infty (a_\ell\cos\ell
x +b_\ell \sin\ell x )\Bigr)^3 {{dx}\over{2\pi}} \cr
&=2\sqrt{2} \int_{\bf T}\cos jx \cos kx\Bigl(
{{a_0}\over{\sqrt{2}}}+\sum_{\ell
=1}^\infty (a_\ell\cos\ell x +b_\ell\sin\ell x
)\Bigr){{dx}\over {2\pi}}&(2.10)\cr}$$
\noindent The matrix that represents 
${{\partial^2 V}\over{\partial
a_j\partial a_k}}$ is 
$$D^2-\beta [v_{jk}]={{7}\over{8}}D^2+{{1}\over{8}}D\Bigl( I-8\beta
\bigl[{{v_{jk}}\over {jk}}\bigr]\Bigr) D\eqno(2.11)$$
\noindent where $D^2\geq I$ and 
$$\eqalignno{\sum_{j,k=1}^\infty &
{{v_{jk}\xi_j\eta_k}\over{jk}}\cr
&=2\sqrt{2}\int_{\bf T}\sum_{j=1}^\infty
{{\xi_j\cos jx}\over{j}}\sum_{k=1}^\infty {{\eta_k\cos
kx}\over{k}}\Bigl({{a_0}\over{\sqrt{2}}}+\sum_{\ell +1}^\infty \bigl(
a_\ell \cos\ell x +b_\ell \sin \ell x\bigr)\Bigr)
{{dx}\over{2\pi}}&(2.12)\cr}$$
\noindent where by the Cauchy--Schwarz inequality
$$\Bigl( \sum_{j=1}^\infty {{\xi_j\cos jx}\over{j}}\Bigr)^2\leq
\sum_{j=1}^\infty {{1}\over{j^2}}\sum_{j=1}^\infty \xi_j^2\leq
{{\pi^2}\over{6}}\sum_{j=1}^\infty \xi_j^2\eqno(2.13)$$
\noindent and 
$$\int_{\bf T}\Bigl({{a_0}\over{\sqrt{2}}}+\sum_{\ell +1}^\infty \bigl(
a_\ell \cos\ell x +b_\ell \sin \ell x\bigr)\Bigr)^2
{{dx}\over{2\pi}}\leq N.\eqno(2.14)$$
\noindent Hence we have
$$ D^2-\beta [v_{jk}]\geq (1/2)D^2;\eqno(2.15)$$
similar estimates apply to the sine terms when we consider 
${{\partial^2 V}\over{\partial
b_j\partial b_k}}$, and to the mixed sine and cosine term which
arise in ${{\partial^2 V}\over{\partial
a_j\partial b_k}}$. Hence $W$ is uniformly convex and thus satisfies Bakry
and Emery's criterion, so $\omega_N^\beta $ satisfies the logarithmic
Sobolev inequality (2.8), and hence $\nu_N^\beta$ satisfies (1.12)
with $\alpha =1/2$.\par
\rightline{$\square$}\par
\vskip.05in

\noindent {\bf Proof of Theorem 1.} We need to extend the
logarithmic Sobolev inequality to a typical pair $N,\beta >0$,
possibly at the expense of a worse constant. So we choose
$K>4\beta \sqrt {N}+1$, and split $\phi\in L^2({\bf T})$ into
the tail $\phi_K(x)=\sum_{k:\vert k\vert \geq K}a_ke^{ikx}$ and the
head 
$h_K(x)=\sum_{k:\vert k\vert \leq K}a_ke^{ikx}$ of the series. We note that for
$\phi\in B_N$, the components satisfy $\Vert h_K\Vert_\infty \leq
(2K+1)^{1/2}N^{1/2}$ and $\int \phi_K(x)^2dx/(2\pi )\leq N$; hence by
some simple estimates 
$$\eqalignno{\Bigl\vert \int_{\bf T}\bigl(\phi (x)^3-\phi_K
(x)^3\bigr){{dx}\over{2\pi}}\Bigr\vert&=\Bigl\vert\int_{\bf
T}\bigl(
3\phi_K(x)^2h_K(x)+3\phi_K(x)h_K(x)^2+h_K(x)^3\bigr)
{{dx}\over{2\pi}}\Bigr\vert&\cr
&\leq 7\beta
(2K+1)^{3/2}N^{3/2}.&(2.16)\cr}$$
\noindent We replace the original potential $V$ by 
$$V_K(a, b)={{a^2_0}\over{2}}+{{1}\over{2}}\sum_{j=1}^\infty j^2(a_j^2+b_j^2)-{{\beta
2\sqrt{2}}\over{6}}\int_{\bf T}\Bigl( \sum_{j=K+1}^\infty 
(a_j\cos jx +b_j\sin jx )\Bigr)^3{{dx}\over {2\pi
}},\eqno(2.17)$$
\noindent which is a bounded perturbation of $V$ on $B_N$ and satisfies 
$$\Vert V-V_K\Vert_\infty \leq 7\beta (2K+1)^{3/2}N^{3/2}.\eqno(2.18)$$
\noindent The matrix $[v_{jk}]_{j,k:\vert j\vert, \vert k\vert\geq
K}$ that arises from $V_K$ via (2.10) involves only high
frequency components and satisfies
$$\beta^2\Bigl\Vert\Bigl[{{v_{jk}}\over{jk}}\Bigr]\Bigr\Vert_{c^2}^2\leq
4\beta^2N\Bigl(\sum_{k=K}^\infty {{1}\over{k^2}}\Bigr)^2
\leq {{4\beta^2N}\over {(K-1)^2}}\leq {{1}\over{4}}\eqno(2.19)$$
\noindent by the choice of $K$. By Lemma 2, $V_K$ is uniformly
convex and the corresponding Gibbs measure satisfies a logarithmic
Sobolev inequality with constant independent of $N$ and $\beta$. \par
\indent Since $V$ is a bounded perturbation of $V_K$, the
Holley--Stroock lemma [7] shows that the Gibbs measure
associated with $V$ also satisfies the logarithmic Sobolev inequality
$$\int_{B_N}F(\phi )^2\log \Bigl(F(\phi
)^2/\int_{B_N}F^2d\nu_N^\beta\Bigr) \nu_N^\beta (d\phi )\leq
4\exp (C\beta^{5/2}N^{9/4})\int_{B_N}\Vert\delta F(\phi )\Vert^2_{L^2}
\nu_N^\beta (d\phi
).\eqno(2.20)$$  
\noindent for come universal constant $C$.\par
\rightline{$\square$}\par

\vskip.05in
\noindent {\bf Corollary 2.} {\sl Let $F:B_N\rightarrow {\bf R}$ be a
Lipschitz function such that $\vert F(\phi )-F(\psi )\vert\leq \Vert
\phi -\psi\Vert_{L^2}$ for all $\phi, \psi\in B_N$; suppose further
that $\int_{B_N}F(\phi )\nu_N^\beta (d\phi )=0$. Then}
$$\int_{B_N}\exp (tF(\phi ))\nu_N^\beta (d\phi )\leq \exp
\Bigl(e^{C\beta^{5/2}N^{9/4}}t^2\Bigr)\qquad (t\in {\bf
R}).\eqno(2.21)$$ 
\vskip.05in
\noindent {\bf Remark.} The result shows that on the probability space $(B_N, \nu_N^\beta )$, the random
variable $F$ has mean zero and takes values that are tightly
concentrated about its mean value.\par
\vskip.05in

\noindent {\bf Proof.} Let $P_n$ be the orthogonal projection onto
${\hbox{span}}\{ e_j:1\leq j\leq n\}$, where $(e_j)$ is some
orthonormal basis of $L^2({\bf T})$. Then $F\circ P_n$ is Lipschitz
continuous on a finite-dimensional subspace, and hence Fr\'echet
differentiable almost everywhere by Rademacher's theorem. We
observe that
$$\Vert\delta F(P_n\phi )\Vert_{L^2}\leq \Vert\nabla F
(P_n\phi )\Vert_{L^2}\leq 1\eqno(2.22)$$
\noindent since $F$ and $P_n$ are Lipschitz. Since
$F\circ P_n\rightarrow F$
uniformly on compact sets as $n\rightarrow\infty$, it suffices by
Fatou's lemma to prove (2.21) for $F\circ
P_n$ and then let $n\rightarrow\infty$.\par
\indent The inequality (2.21) then follows from Theorem 1 by the general theory of
functional inequalities, as in [16]; here we give a brief argument. Let
$$J(t)=\int_{B_N}e^{tF(\phi )}\nu_N^\beta (d\phi )\eqno(2.23)$$
\noindent which defines an analytic function of $t$ such that $J(0)=1$
and $J'(0)=0.$ Further, the logarithmic Sobolev inequality gives
$$\eqalignno{ tJ'(t)&=
\int_{B_N}tF(\phi )e^{tF(\phi )}\nu_N^\beta (d\phi )\cr
&\leq J(t)\log J(t) +{{t^2}\over{2\alpha }}J(t)\qquad
(t>0),&(2.24)\cr}$$
\noindent which integrates to the inequality
$$J(t)\leq \exp \Bigl( {{t^2}\over{2\alpha }}\Bigr).\eqno(2.25)$$
\rightline{$\square$}\par
\vskip.05in

\noindent {\bf Remark.} We leave it as an open problem to determine whether
$\nu_N^\beta$ satisfies (1.12) with a constant independent of $N$ for
given $\beta$.\par

\vskip.1in
\noindent {\bf 3. Stationary points of the Hamiltonian on spheres}\par
\vskip.05in
\noindent The Hamiltonian $H(\phi )$ is unbounded above and below for
$\phi\in L^2({\bf T})$; however, we can consider the minimal energy
constrained to the spheres in $L^2({\bf T})$:
$$E_N=\inf\Bigl\{ H(\phi) :\int_{\bf T}\phi (x)^2{{dx
}\over{2\pi}}=N\Bigr\}.\eqno(3.1)$$
\indent Korteweg and de Vries introduced a travelling
wave solution $u(x,t)=v(x-ct)$ of $u_t+u_{xxx}+\beta uu_x=0$ which
is periodic and is commonly known as the cnoidal wave. We
recover this solution below. Subject 
to some reservations, Drazin [5,6] showed that the
cnoidal wave is linearly stable with respect to any infinitesimal perturbation.\par
We recall from [13] Jacobi's sinus amplitudinus of modulus $k$ is 
${\hbox{sn}}(x\mid k)=\sin \psi$ where
$$x=\int_0^\psi {{d\theta}\over{\sqrt{1-k^2\sin^2\theta}}}.\eqno(3.2)$$
\indent For $0<k<1$, let $K (k)$ be the complete elliptic
integral
$$K(k)=\int_0^{\pi/2}{{dt}\over{\sqrt{1-k^2\sin^2t}}};\eqno(3.3)$$
\noindent then ${\hbox{sn}}(z\mid k)^2$
has real period $K$ and complex period $2iK(\sqrt{1-k^2})$.\par
\indent For $\ell >0$, the standard form of
Lam\'e's equation in [11] is
$$\Bigl( -{{d^2}\over{dz^2}}+\ell (\ell +1)k^2{\hbox{sn}}(z\mid
k)^2\Bigr)\Phi (z)=\mu \Phi (z).\eqno(3.4)$$
\noindent The $L^2$ spectrum of (3.4) is determined by a sequence
$\lambda_0<\lambda_1\leq \lambda_2<\lambda_3\leq \lambda_4<\dots, $
which is infinite except for $\ell =1, 2, \dots$. Typically 
$\sigma_B=\cup_{j=0}^\infty
[\lambda_{2j}, \lambda_{2j+1}]$ gives the Bloch spectrum, so that for
$\mu\in \sigma_B$ there exists a bounded solution to (3.4);
whereas for $\mu\in (-\infty , \lambda_0)\cup\bigcup_{j=0}^\infty (\lambda_{2j+1},
\lambda_{2j+2})$ all nontrivial solutions of (3.4) are unbounded, and
we say that $\mu$ belongs to an interval of instability. In the special case of $\ell =1, 2,\dots, $ there
are only $\ell +1$ intervals of instability, namely $(\lambda_{2j+1},
\lambda_{2j+2})$ $(j=0, \dots ,j-1)$ and $(-\infty ,
\lambda_0)$.\par

\vskip.05in

\noindent {\bf Theorem 3.} {\sl (i) Let 
$\phi \in C^2({\bf T};{\bf R})$ be a stationary point for the energy
$$H_\lambda (\phi ) ={{1}\over{2}}\int_{\bf T} \phi'(x)^2{{dx}\over{2\pi}}-{{\beta}\over{6}}\int_{\bf T} \phi
(x)^3{{dx }\over{2\pi}}-{{\lambda }\over{2}}\int_{\bf T}
\phi (x)^2{{dx}\over{2\pi}}.\eqno(3.5)$$
Then $\phi$ satisfies the differential equation
$$\phi''(x)+{{\beta}\over{2}}\phi (x)^2 +\lambda \phi (x)=0,\eqno(3.6)$$
\noindent so either $\phi$ is constant or an elliptic function. \par
\indent (ii) For $\beta>0$, the energy $H_\lambda$ on $S_N$ has a local minimum at 
$\phi =-\sqrt{N}$.\par
\indent (iii) Let $\phi$ be the elliptic function} 
 $$\phi(x)=f_1-(f_1-f_2)\Bigl[{\hbox{sn}}\Bigl(\sqrt{{{\beta (f_1-f_3)}\over
{12}}}(x_1-x)\Bigl\vert\sqrt{{{f_1-f_2}\over{f_1-f_3}}}\Bigr)\Bigr]^2\eqno(3.7)$$
\noindent {\sl for suitable real constants $f_3<f_2<f_1$ and
$\phi (x_1)=f_1$. Let $(-\infty ,\lambda_0]$ be the zeroth order
interval of instability of Lam\'e's equation 
$$y''(x)+{{\beta}}\phi(x)y(x)+\lambda y(x)=0.\eqno(3.8)$$
\noindent Then for $\lambda <\lambda_0$ the energy $H_\lambda $ has a 
local minimum at $\phi$; whereas for $\lambda >\lambda_0$ the
stationary point is neither a local maximum nor a local minimum.}\par

\vskip.05in
\noindent {\bf Proof.} (i) We suppose
that $\beta >0$. One can easily expand $H_\lambda (\phi +t\psi )$ as a
cubic polynomial in $t$ and examine the conditions that ensure that
$t=0$ gives a local minimum. The equation ${{\delta H_\lambda }\over{\delta\phi}}=0$
reduces to the differential equation (3.6) which has constant solutions 
$\phi =0$ (which does not belong to $S_N$) and $\phi = -2\lambda
/\beta$, and a non-constant solution satisfying
$${{1}\over{2}}\phi'(x)^2 +{{\beta}\over{6}}\phi
(x)^3+{{\lambda}\over{2}}\phi (x)^2=C\eqno(3.9)$$
\noindent with $C$ some constant. This equation has periodic solutions if
and only if $\beta^2
<3\lambda^3/(2C)$; equivalently, for such constants there exist real roots  $f_3<f_2<f_1$ such that
$$-{{\beta}\over{6}}\phi^3-{{\lambda}\over{2}}\phi^2+C
=-{{\beta }\over{6}}(\phi -f_1)(\phi
-f_2)(\phi -f_3).\eqno(3.10)$$ 
\noindent To find these roots, we introduce $z=1/\phi$, which
satisfies the cubic $z^3-\lambda z/(2C)-\beta /(6C)=0$ with
discriminant
$$D={{-\lambda^3}\over {216 C^3}}+{{\beta^2}\over {144
C^2}};\eqno(3.11)$$
\noindent so by writing $re^{i\theta} =\beta
/(12C)+i\sqrt {-D}$, we have 
$${{1}\over{f_1}}=2r^{1/3}\cos {{\theta}\over{3}},\quad  
{{1}\over{f_2}}=2r^{1/3}\cos {{\theta
+2\pi}\over{3}}, \quad {{1}\over{f_3}}=2r^{1/3}\cos {{\theta
+4\pi}\over{3}},\eqno(3.12)$$
\noindent for some choice of the polar angle. To convert to the
standard form (3.4) of Lam\'e's equation, we introduce
$$k=\sqrt{{{f_1-f_2}\over{f_1-f_3}}}, \qquad k^2\ell (\ell +1) =\beta
(f_1-f_3),\quad \gamma =\cos {{\theta }\over{3}},\eqno(3.13) $$
\noindent where $0<k<1$ and $\ell >0$, and by some 
trigonometry deduce that
$$k^2={{2\gamma\sqrt {1-\gamma^2}}\over
{\sqrt{3}/2-\sqrt{3}\gamma^2-\gamma\sqrt{1-\gamma^2}}},\eqno(3.14)$$
\noindent which is an algebraic function of the parameters. Using
the definition (3.2), one can show that the 
solution of (3.9) is given by
 $$\phi (x)=f_1-(f_1-f_2)\Bigl[{\hbox{sn}}\Bigl(\sqrt{{{\beta (f_1-f_3)}\over
{12}}}(x_1-x)\Bigl\vert k
\Bigr)\Bigr]^2;\eqno(3.15)$$
\noindent where $\phi$ is $2\pi$-periodic provided that 
$$2\pi \sqrt{{{\beta (f_1-f_3)}\over
{12}}}=2K (k).\eqno(3.16)$$ 
\indent (ii) For any stationary point $\phi$,
we have
$$H_\lambda (\phi +t\psi )=H_\lambda (\phi )+{{t^2}\over{2}}\int_{\bf T}\bigl(
\psi'(x)^2-(\lambda +\beta \phi (x))\psi
(x)^2\bigr){{dx}\over{2\pi}}-{{\beta t^3}\over{6}}\int_{\bf T}\psi
(x)^3{{dx}\over{2\pi}}.\eqno(3.17)$$
\indent First we deal with the constant stationary points. Note
that
$\phi =-\sqrt{N}$ gives a local minimum for $H_\lambda $ with
$\lambda =\beta\sqrt{N}/2$ and $H_\lambda (-\sqrt {N})=-\beta N^{3/2}/12
$ and the term of order $t^2$ is positive. The other constant solution $\phi =\sqrt{N}$ has $\lambda =-\beta
\sqrt{N}$ and $H_\lambda (\sqrt N)=\beta N^{3/2}/12.$\par
\indent (iii) Next we take $\phi$ from
(3.9), so that 
$$H_\lambda (\phi )=-{{\beta}\over{3}}\int_0^{2\pi}\phi (t
)^3{{dt}\over{2\pi}}-\lambda\int_0^{2\pi} \phi (t)^2
{{dt}\over{2\pi}}+C.\eqno(3.18)$$
\noindent To determine when this is a local minimum, we consider 
a nontrivial solution of the equation 
$$y''(x)+{{\beta}}\phi (x)y(x)+\lambda y(x)=0,\eqno(3.19)$$
\noindent which is the linearization of (3.6), and we recognise this
as a form of Lam\'e's equation. \par
\indent First
suppose that $y$ has only finitely many zeros on ${\bf R}$;
then by Hamel's theorem [10] 
$$\int_{\bf T}\bigl(
\psi'(x)^2-(\lambda +\beta \phi (x))\psi
(x)^2\bigr){{dx}\over{2\pi}} \geq 0\eqno(3.20)$$
\noindent holds for all $2\pi$ periodic and continuously differentiable
functions $\psi$. Let $\lambda_0$ be the supremum of such $\lambda$. For $\lambda <\lambda_0$ we have
strict inequality for in (3.20) all nonzero $\psi$, so $H_\lambda $ has a local
minimum at $\phi$.\par

\indent Suppose contrariwise that $y$ has infinitely many zeros; so by Hamel's
theorem there exists $\psi$ such that    
$$\int_{\bf T}\bigl(
\psi'(x)^2-(\lambda +\beta \phi (x))\psi
(x)^2\bigr){{dx}\over{2\pi}} <0;\eqno(3.21)$$
\noindent hence $H_\lambda $ does not have a local minimum at $\phi.$ In
particular, this
happens when $\lambda >\lambda_0$.\par
\rightline{$\square$}\par
\vskip.1in
\noindent {\bf 4. The modified periodic KdV equation}\par
\vskip.05in

\noindent Lebowitz, Rose and Speer [8,9] considered Hamiltonians
$$H_p(\phi )={{1}\over{2}}\int_{\bf T}\phi'
(x)^2{{dx}\over{2\pi}}-{{\beta}\over{p(p-1)}}\int_{\bf T}\vert \phi
(x)\vert^p{{dx}\over{2\pi}}\eqno(4.1)$$
\noindent and showed that there exists a Gibbs measure
with potential $H_p(\phi )$ on $\{ \phi\in L^2({\bf T}): \Vert
\phi\Vert^2_{L^2}\leq N\}$ as in (1.6) for $\beta , N>0$ and
$2\leq p<6$, but
not for $p>6.$ When $\beta <0$, the Hamiltonian is non focusing
and the problem of normalizing the probability measures is easy to
address.\par 
\indent With this Hamiltonian and $\beta >0$, they introduced the Gibbs
measure
$$\nu_N^\beta (d\phi )=Z_N(\beta )^{-1}{\bf I}_{B_N}(\phi )e^{-H_4(\phi
)}\prod_{e^{ix}\in {\bf T}}d\phi (x),\eqno(4.2)$$
\noindent and called it the modified canonical ensemble.\par
\indent Their analysis of the case $p=4$ showed that for small $\beta
>0$, the constant solution was stable; whereas for large $\beta >0$,
the soliton solution was stable. In this section we strengthen
their result by showing that for $p=4$ the modified canonical ensemble
satisfies a logarithmic Sobolev inequality when $\beta N$ is small. The case $p=4$ corresponds to the modified KdV equation
$${{\partial \phi}\over{\partial
t}}+{{\partial^3\phi}\over{\partial\theta^3}}+\beta\phi^2{{\partial\phi}
\over{\partial\theta}}=0;\eqno(4.3)$$
\noindent so the result suggests that at low temperatures, solutions of
the mKdV equation are most likely to occur near to the ground state.\par

\vskip.05in
\noindent {\bf Theorem 4.} {\sl There exists $C>0$ such that for 
$0<\beta N<C$ the Gibbs measure $\nu_N^\beta$ satisfies the logarithmic Sobolev inequality with
 $\alpha =1/2$.}
\vskip.05in
\noindent {\bf Proof.} This is similar to the proof of Lemma 2. To simplify notation, we consider the
Hamiltonian 
$$H(a_n)=\sum_{n=1}^\infty n^2a_n^2-{{\beta}\over{12}}\int_{\bf T}
\Bigl( \sum_{n=1}^\infty a_n\cos n\theta \Bigr)^4
{{d\theta}\over{2\pi}}\eqno(4.5)$$
\noindent on $B_N=\{(a_n)\in \ell^2({\bf R}) : \sum_{n=1}^\infty
a_n^2\leq N\}$ which has essentially the same properties as the true
Hamiltonian in the Fourier components. The corresponding Hessian matrix has entries
$${{\partial^2H}\over{\partial a_j\partial a_\ell}}=\delta_{j\ell}
j^2-\beta\int_{\bf T}\Bigl(\sum_{n=1}^\infty a_n\cos n\theta \Bigr)^2
\cos j\theta \cos \ell \theta {{d\theta}\over{2\pi}},\eqno(4.6)$$
\noindent and we deduce that 
$$\sum_{j,\ell =1}^\infty {{\partial^2H}\over{\partial a_j\partial
a_\ell}}{{\xi_j\xi_\ell}\over{j\ell}}=\sum_{j=1}^\infty\xi_j^2 -\beta
\int_{\bf T}\Bigl(\sum_{j=1}^\infty a_n\cos n\theta \Bigr)^2
\sum_{j=1}^\infty {{\xi_j\cos j\theta}\over{j}}\sum_{\ell
=1}^\infty {{\xi_\ell \cos\ell\theta}\over{\ell}}
{{d\theta}\over{2\pi}},\eqno(4.7)$$
\noindent where by the Cauchy--Schwarz inequality 
$$\Bigl(\sum_{j=1}^\infty {{\xi_j\cos j\theta}\over{j}}\Bigr)^2\leq 
\sum_{j=1}^\infty {{1}\over{j^2}} \sum_{j=1}^\infty
\xi_j^2;\eqno(4.8)$$
\noindent hence 
$$\sum_{j,\ell =1}^\infty {{\partial^2H}\over{\partial a_j\partial
a_\ell}}{{\xi_j\xi_\ell}\over{j\ell}}\geq
{{1}\over{2}}\sum_{j=1}^\infty \xi_j^2\eqno(4.9)$$
\noindent and Bakry and Emery's condition is satisfied. The logarithmic
Sobolev inequality follows.\par
\rightline{$\square$}\par
\vskip.05in
\noindent {\bf Remark.} The author has not succeeded in proving a logarithmic Sobolev
inequality for $p=4$ and $N\beta$ large due to the lack of a suitable substitute for
(2.16).\par

\vskip.1in

\noindent {\bf 5. Periodic solutions of the cubic nonlinear Schr\"odinger
equation}\par
\vskip.05in
    
\indent The Hamiltonian 
$$H={{1}\over{2}}\int_{\bf T}\bigl(
(Q')^2+(P')^2\bigr)\, dx +{{\beta}\over{4}}\int_{\bf
T}\bigl(Q^2+P^2\bigr)^2\,dx\eqno(5.1)$$
\noindent gives the canonical equations of motion
$$\left[\matrix{0&-1\cr 1&0\cr}\right] {{\partial }\over{\partial t}}\left[\matrix {Q\cr P\cr}\right] 
=-{{\partial^2 }\over{\partial x^2}}\left[\matrix {Q\cr P\cr}\right]
+\beta (Q^2+P^2)\left[\matrix {Q\cr P\cr}\right]\eqno(5.2)$$
\noindent which give the cubic Schr\"odinger equation for $u=P+iQ$ 
$$-iu_t=-u_{xx}+\beta \vert u\vert^2u\eqno(5.3)$$
\noindent over the circle. Again we take $\beta >0$, which gives the
focusing case. The appropriate number operator is
represented by $N={{1}\over{2}}\int_{\bf T}(P^2+Q^2) dx,$ which is 
invariant under the flow, and we introduce 
$$B_N=\{
u=P+iQ:\int_{\bf T}( P^2+Q^2)dx\leq N\}.\eqno(5.4)$$
\noindent As discussed in [15], one can normalize the
measure
$$\nu_N^\beta (du)=Z_N(\beta )^{-1}\exp\bigl( -H(u)\bigr){\bf I}_{B_N}(u)\prod_{x\in {\bf T}}
du(x)\eqno(5.5)$$  
\noindent so that it gives a probability measure on $B_N$. Bourgain
has shown that $\nu_N^\beta$ is invariant under the flow associated
with the cubic Schr\"odinger equation; see [4, p.124] \par
\vskip.05in
\noindent {\bf Theorem 5.} {\sl There exists $C>0$ such that
for all $0<\beta N<1/2$ the Gibbs measure $\nu_N^\beta$ satisfies the logarithmic Sobolev inequality with
 $\alpha =1/2$.} 
\vskip.05in
\noindent {\bf Proof.} This follows the proof of Lemma 2 and Theorem 4
closely, hence is omitted.\par 
\rightline{$\square$}\par

\vskip.05in
\noindent {\bf References}\par
\noindent [1] D. Bakry and M. Emery, Diffusions hypercontractives.
Seminaire de probabilit\'es XIX 1983/84 177--206 Lecture Notes in
Mathematics, (Springer, Berlin, 1985).\par
\noindent [2] J. Bourgain, Fourier transform restriction phenomena for
certain lattice subsets and applications to nonlinear evolution
equations II. The KdV equation. {Geom. Funct. Anal.} {3}
(1993), 209--262.\par
\noindent [3] J. Bourgain, Periodic nonlinear Schr\"odinger equation and
invariant measures, {Commun. Math. Phys.} {166} (1994),
1--26.\par
\noindent [4] J. Bourgain, {Global solutions of nonlinear
Schr\"odinger equations}, (American Mathematical Society, 1999).\par
\noindent [5] P.G. Drazin, On the stability of cnoidal waves, {Quart.
J. Mech. Appl. Math.} {30} (1977), 91--105.\par
\noindent [6] P.G. Drazin, {Solitons}, (Cambridge University Press,
1983).\par
\noindent [7] R. Holley and D. Stroock, Logarithmic Sobolev inequalities
and stochastic Ising models, {J. Statist. Phys.} {46}, 1987,
1159--1194.\par 
\noindent [8] J.L. Lebowitz, H.A. Rose and E.R. Speer, Statistical
mechanics of the nonlinear\par
\noindent Schr\"odinger equation, {J. Statist.
Phys.} {50} (1988), 657--687.\par
\noindent [9] J.L. Lebowitz, H.A. Rose and E.R. Speer, Statistical
mechanics of the nonlinear\par
\noindent Schr\"odinger equation II: mean field approximation, {J. Statist.
Phys.} {54} (1989), 17--56.\par
\noindent [10] W. Magnus and S. Winkler, {Hill's equation}, (Dover
Publications, New York, 1966).\par
\noindent [11] R.S. Maier, Lam\'e polynomials, hyperelliptic
reductions and Lam\'e band structure, {Phil. Trans. R. Soc. A}
{336} (2008), 1115--1153.\par
 \noindent [12] H.P. McKean, Geometry of differential space, {Ann. 
Probability} {1} (1973), 197--206.\par
\noindent [13] H. McKean and V. Moll, {\sl Elliptic curves: function
theory, geometry, arithmetic,} (Cambridge University Press, 1997).\par
\noindent [14] H.P. McKean and K.L. Vaninsky, Brownian motion with restoring
drift: the petit and microcanonical ensembles, {Comm. Math. Phys.}
{160} (1994), 615--630.\par
\noindent [15] H.P. McKean and K.L. Vaninsky, Action-angle variables
for the cubic nonlinear\par
\noindent Schr\"odinger equation, {Comm. Pure Appl.
Math.} {50} (1997), 489--562.\par
\noindent [16] C. Villani, {Topics in Optimal Transportation,} (American
Mathematical Society, 2003).\par

\vfill
\eject
\end